\input amstex
\documentstyle{amsppt}
\redefine\c{\hat\Bbb C}
\redefine\C{\hat\Bbb C}
\redefine\D{\Bbb D}
\define\Ch{{\Cal C}}

\define\ce{Collet-Eckmann}
\define\cec{Collet-Eckmann condition}
\define\fmn{f^{-n}}
\redefine\R{\Bbb R}
\redefine\N{\Bbb N}
\redefine\B{{\Cal B}}
\define\T{{\Cal T}}

\define\a{\alpha}
\redefine\b{\beta}
\redefine\d{\delta}
\redefine\l{\lambda}

\define\e{\varepsilon}
\define\m{\mu}
\define\diam{\operatorname{diam}}              
\define\dist{\operatorname{dist}}
\define\karo{\hfill$\square$}

\NoBlackBoxes

\document
\topmatter
\title Porosity of Collet-Eckmann Julia sets
\endtitle
\author Feliks Przytycki and Steffen Rohde 
\endauthor

\date   April 18, 1996 \enddate

\address{Institute of Mathematics, Polish Academy of Sciences, 
Sniadeckich 8, 00-950 Warszawa, Poland}
\endaddress
\address{Department of Mathematics, TU Berlin,
Stra\ss e d. 17. Juni 136, 10623 Berlin, Germany}
\endaddress

\abstract We prove that the Julia set of a rational map of the Riemann
sphere satisfying the Collet-Eckmann condition 
and having no parabolic periodic point is mean porous, if it is
not the whole sphere. It follows that the Minkowski dimension of the
Julia set is less than 2.
\endabstract

\thanks The authors acknowledge support by Polish KBN 
Grant 2 P301 01307.
\endthanks

\endtopmatter


\heading1. Introduction \endheading
Let $f:\c\to\c$ be a rational map. Then $f$ is said to satisfy the {\it Collet-
Eckmann condition} if there are constants $C>0$ and $\l>1$ such that
$$|(f^n)'(f(c))|\geq C \l^n\tag CE$$
for all $n$ and
all critical points $c\in J(f)$ of $f$ whose forward orbit does not meet another
critical point ($J(f)$ stands for the Julia set of $f$). 
Here and in what follows derivatives and distances are always with
respect to the spherical metric of $\C,$ unless stated otherwise.

A set $E\subset\c$ is called {\it mean porous} if there are constants
$p_1<\infty$ and $p_2>0$ such that for each $z\in E$ the following holds:
There is an increasing sequence $n_j$ of integers and points
$z_j$ with $\dist (z,z_j) \leq 2^{-n_j}$ such that $n_j<p_1 j$ and
$\dist(z_j,E)>p_2 2^{-n_j}.$
Roughly speaking, the scales in which $E^c$ contains a disc of  size
proportional to the scale
have a density bounded uniformly from below.

\proclaim{Theorem 1.1} If $f$ satisfies the \cec , has no parabolic
periodic point and if $J(f)\neq\c,$ then $J(f)$ is mean porous.
\endproclaim

In \cite{KR} it was proved that mean porous sets on the sphere
have Minkowski dimension $<2,$ see also Section 4.
As an immediate consequence we obtain

\proclaim{Corollary 1.2} Under the assumptions of Theorem 1.1,
the Minkowski dimension of $J(f)$ is less than 2.
\endproclaim
 
Maps satisfying (CE) were first considered by Collet and Eckmann in \cite{CE}.
Benedicks and Carleson showed in \cite{BC, Theorem 1} 
that the set of (real) parameters
$c$ for which $z^2+c$ satisfies (CE) is of positive measure. Nowicki and the 
first author showed in \cite{NP}
H\"older conjugacy with tent map for Collet-Eckmann
interval maps and conjectured that the basin of $\infty$ is a H\"older domain
(for complex quadratic maps with (CE)). A new tool to control distortion
of components of preimages of discs, {\it shrinking neighborhoods} in the terminology
of Graczyk and Smirnov, was introduced by the first author
in \cite{P1} (see also \cite{PUZ, p. 198}). 
There it was used
to establish (for $f$ as in Theorem 1.1, say)
existence of absolutely continuous invariant measures for
conformal measures. It was also shown that Hausdorff dimension, Minkowski
dimension and hyperbolic dimension of the Julia set coincide.
In \cite{P2}, a question of Bishop and Lyubich was answered negatively
by showing that the Hausdorff dimension of maps as in Theorem 1.1 is
less than 2, provided some additional condition (M. Tsujii condition) holds.
Shrinking neighborhoods were used to go from small scale to large scale
for $f(c),$ $c$ critical.
Other types of counterexamples come from Graczyk's work on real
quadratic Fibonacci polynomials \cite{G}
and from McMullen's work \cite{McM}.

In the summer 1995 the first author had a discussion with M. Lyubich who
suggested trying to estimate the dimension of the Julia set of Collet-Eckmann
maps directly (without going through conformal measure as in \cite{P2}), as
going from large scale to small scale should 'pull back a hole' (i.e. a disc
contained in $J^c$ on large scale) to small scales. It was realized that
some notion of porosity could be involved. 

Koskela and the second author introduced in \cite{KR} the notion of 
mean porosity, showed that
the boundary of a H\"older domain is mean porous and proved 
(based on work of Jones and Makarov \cite{JM}
and of Smith and Stegenga \cite{SS}) that the Minkowski
dimension of mean porous sets is $<2.$ 

Graczyk and Smirnov proved in \cite{GS}
that the components of the Fatou set of \ce\ maps are H\"older domains.
They concluded that, for polynomial \ce\ maps, the Minkowski dimension of the
Julia set is $<2$ and asked whether it is always $<2.$
Nazarov, Popovici and Volberg \cite{NPV} extended the
results of a first version of \cite{GS} to disconnected Julia sets 
of polynomials and
raised the question whether the Julia set could be always mean porous.

In the present paper we give a positive answer to this question.
We carry out the program of 'pulling back holes from
large to small scale'. The main ingredient, besides shrinking neighborhoods,
is an estimate from \cite{DPU} on the average distance of an orbit $f^n(x)$
from the set of critical points of $f,$ which is our substitute for the
Tsujii condition used in \cite{P2}.
As a byproduct of our proof of mean porosity, we obtain a new approach 
to the Graczyk and Smirnov H\"older theorem, outlined in Section 3.

Pulling back holes is done with uniformly 
bounded criticality on surrounding discs, so
in an abundance of scales \ce\ maps behave like semihyperbolic maps
\cite{CJY}\cite{DU}\cite{U}.

In Section 4 we define a notion of porosity that is slightly stronger than
the above mean porosity, show that Collet-Eckmann Julia sets ($\neq\c,$ without
parabolic periodic points) satisfy this condition and give a simple proof that
the Minkowski dimension of such sets is less than 2.
 

\heading2. Proof of Theorem 1\endheading

Consider a disc $B=B(x,\d)$ of spherical radius $\d$ around $x\in\C$
and a connected component $W$ of 
$f^{-n}(B).$ We are mainly interested in the case that $f^n_{|W}$ has
at most $D$ critical points (counted with multiplicity), where $D$ is some
fixed number. In this situation, say that $f^n$ is {\it D-critical on W}.
Then $f^n_{|W}$ has distortion properties similar to those of conformal maps. 
We collect the  estimates needed in our paper as Lemma 2.1 below
(see \cite{P1, Section 1}).  The radius
$\d$ will always be assumed to be less than $\diam(\C)/2.$

\proclaim{Lemma 2.1} For each $\e>0$ and $D<\infty$ there  are constants
$C_1$ and $C_2$, such that the following holds for all rational maps
$F:\C\to\C,$ all $x\in\C$ and all $t$ with $1/2\leq t<1:$

Assume that $W$ resp. $W'$ are simply connected components of 
$F^{-1}(B(x,\d))$ resp. $F^{-1}(B(x,t\d))$ with $W\supset W'$.
Assume further that $\C\setminus W$ contains a disc of radius $\e$
and that $F$ is D-critical on $W.$ Then
$$|F'(y)| \diam(W') \leq C_1 (1-t)^{-C_2}\d\tag2.1$$
for all $y\in W'.$

Furthermore, if $t=1/2$ and $0<\tau<1/2$,
let $B''=B(z,\tau\d)$ be any disc contained in $B(x,\d/2)(=F(W'))$
and let $W''$ be a component of $F^{-1}(B'')$ contained in $W'.$ Then
$$\diam(W'')\leq C_3\diam(W')\tag 2.2$$
with $C_3=C_3(\tau,\e,D)$ and $C_3\to 0$ as $\tau\to 0$ (for fixed $\e,D$).
Finally,
$$ W'' \text{ contains a disc of radius } \geq C_4\diam(W')\tag2.3$$
around every preimage of $F^{-1}(z)$ that is contained in $W''.$
Here $C_4=C_4(\tau,\e,D).$
\endproclaim

\demo{Proof}
We will give a short proof of (2.1) (which is essentially Lemma 1.4 in
\cite{P1}; the statement in the preprint version of \cite{P1} is
imprecise) and of (2.3). The inequality (2.2) is 
\cite{P1, (1.5') with $\l=1/2$}.

We may assume that the disc of radius $\e$ contained
in $\C\setminus W$ is centered at $\infty$, hence $W$ is a simply
connected planar domain bounded by some constant depending on $\e$ only.
We may further assume that $B(x,\d)$ is the unit disc $\D$ and that
$B'=B(x,t \d)$ is the disc $\{|w|<t\}.$ Finally we assume (translate
$W'$ if necessary) $0\in W'$ and $F(0)=0.$ 

To prove (2.1), we need to show that
$$|F'(y)|\diam(W')\leq C_1 (1-t)^{-C_2}$$
for each $y\in W',$
where now (and during the rest of the proof of Lemma 2.1) diameters and 
derivatives are with respect to the {\it euclidean} metric
(here is where we need the assumption involving $\e$).

Let $g:\D\to W$ be a conformal map with $g(0)=0$ and set $h=F\circ g.$
Then $h$ is a Blaschke product of degree $\leq D,$ and $h(0)=0.$
We will show that 
$$G:=g^{-1}(W')\subset \{|w|<1-C_1'(1-t)^{C_2'}\}$$ 
for constants $C_1',C_2'$ depending only on $D.$ 
>From this (2.1) follows
immediately since (writing $d=1-C_1'(1-t)^{C_2'}$ for short)
$\diam(W')\leq 2\frac{|g'(0)|}{(1-d)^2}$
by Koebe distortion \cite{P,Chapter 1.3},
and (for $y\in W'$ and $u\in G$ with $g(u)=y$)

$|F'(y)|=|h'(u)|/|g'(u)|
  <8[1/(1-|u|^2)]/[(1-|u|)|g'(0)|],$
\newline\noindent
using Koebe distortion again.

To prove $G\subset \{|w|<d\},$ write 
$h(u)= \prod_{n=1}^{deg(h)} [(u-a_n)/(1-\overline{a_n}u)]
= \prod_{n=1}^{deg(h)}T_n(u)$ and
notice that, if $|h(u)|<t,$ then at least one of the factors has to
be of absolute value
$< t^{1/D}.$ Denoting the hyperbolic metric of $\D$ by $\rho$ and using
$\rho(u,a_n)=\log((1+|T_n(u)|)/(1-|T_n(u)|)$, easy calculation shows
$\rho(u,a_n)< C(D)+\log(1/(1-t))$ for such $n.$ In other words, every
$u\in G$ has hyperbolic distance $\leq C(D)+\log(1/(1-t))$ from the set
$\{a_1,...,a_{deg(h)}\}.$ Since $0\in G$ and $deg(h)\leq D,$ it easily follows
that $G\subset \{w\in \D:\rho(0,w)<2D(C(D)+\log(1/(1-t)))\}
\subset \{|w|<1-C_1'(1-t)^{C_2'}\}$ and (2.1) is proven.

Now (2.3) is a simple consequence of the above and Koebe distortion:
Using notation as above, we already know $G\subset\{|w|<d\},$ and by the
Lemma of Schwarz $G$ contains a hyperbolic disc around every preimage
$h^{-1}(z)\in G,$ of hyperbolic radius at least the hyperbolic radius of $B''.$
\karo\enddemo

Consider a rational map $f$ with $J(f)\neq\C$ and without parabolic periodic points.

>From now on, we will always assume that $\d$ small enough to guarantee
that all components of $f^{-n}(B(y,\d))$ do not meet a fixed open set
that contains all critical points of $f$ in the Fatou set,
for all $n$ and all $y\in J$ (pick any small open neighborhood 
of the critical points in the Fatou set
and take the union of its forward orbit under $f$).
This allows us to apply Lemma 2.1 to simply connected $D-$ critical components,
with bounds depending only on $D.$

Now fix $\d>0$ and $D<\infty$ and consider
$B=B(f^n(x),\d)$ together with the component $W$ of $\fmn(B)$ containing $x.$
We call $n$ a {\it good time} for $x$ and denote the set of good times
by $G(x),$ if $f^n$ is $D-$ critical on $W.$

\proclaim{Lemma 2.2 (uniform density of good times)} There exists $\d>0$ and
$D<\infty$ such that the lower density of $G(x)$ in $\N$ is at least
$1/2,$
$$\inf_{n} \frac{\#(G(x)\cap[1,n])}{n} \geq \frac12.$$
\endproclaim

\demo{Proof} The proof is a modification of the proof of Lemma 2.1
in \cite{P2}. There the Tsujii condition was used to obtain times $n$ where
$f^n:W\to B(f^n(c),\d)$ has degree one. The main ingredient here is
the inequality (3.3) of \cite{DPU}, an estimate of the average distance from
critical points.

Fix $x\in J.$ As in \cite{P2}, set 
$$\phi(n)=-\log(\dist(f^n(x), Crit(f,J)),$$
where $Crit(f,J)$ denotes the set of critical points of $f$ that are
contained in $J.$ Then by (3.3) of \cite{DPU} there exists a constant
$C_f,$ such that for each $n\geq1$
$${\sum_{j=0}^n}\ ' \phi(j) \leq n C_f ,\tag 2.4$$
where $\sum '$ denotes summation over all but at most $\#Crit(f,J)$
indices.

One could view the 'graph' of $\phi$ as the union of all vertical
line segments $\{n\}\times [0,\phi(n)]$ in $\Bbb R^2.$ Then each segment
throws a {\it shadow} $S_n = (n, n+\phi(n) K_f] \subset \R,$
where we set $K_f=2\nu /\log(\l)$ and denote by $\nu$ the largest degree
of all critical points in $J$ of all iterates of $f.$

The shadows of the exceptional indices in (2.4) could be infinitely long,
but nevertheless (2.4) implies that many of the times $n$ belong to boundedly
many shadows: Indeed, set $N_f=2(\#Crit(f,J)+ C_f K_f)$ and
$$A=\{j\in\N : j \text{ belongs to at most } N_f \text{ shadows } S_n\},$$
then for each $n$ we obtain from (2.4)
$${\sum_{j=0}^n}\ ' |S_j| \leq C_f K_f n$$
and conclude that
$$\frac{\#(A\cap [0,n])}{n}\geq\frac12.$$
We now show that each $n\in A$ is a good time for $x$, i.e.
$A\subset G(x)$ (with $D=\nu N_f$ and $\d$ suitable).
We use the technique of \cite{P1} of
'shrinking neighborhoods'.

Fix once and for all a subexponentially decreasing sequence $b_j>0$
with $\prod_{j=1}^\infty (1-b_j) > 1/2.$ Fix $n\in A$ and consider
the sequence 
$$B_s=B(f^n(x),2\d\prod_{j=1}^s (1-b_j))$$
of neighborhoods of $B=B(f^n(x),\d),$ together with (compatible)
connected components $W_s$ of $f^{-s}(B_s)$ and
$W_s'$ of $f^{-s}(B_{s+1}).$

Recall the main idea of shrinking neighborhoods from \cite{P1}:
If (along backwards iteration from $f^n(x)$) a critical value is met
in $W_s$ but not in $W_s',$ then it can be ignored (because
$f^{t-s}$ maps $W_t$ into $W_s'$ for $t>s$). As $W_s'$ sits
'well inside' $W_s,$ distortion on $W_s'$ can be controlled.

We want to show that if $W_s$ contains a critical point, then
$n$ belongs to the shadow $S_{n-s}.$ Assume this is not the case.
Then there is a smallest such s, a critical point $c\in W_s,$ and
$f^{s-1}:W_{s-1}\to B_{s-1}$ is at most $\nu N_f$-critical (as $n\in A$
and $s$ is smallest).  
It is not hard to see that $W_{s-1}$ is simply connected.
Use induction: 
The fact that at most $N_f$ of the domains $W_{s-2},...,W_1$ 
contain critical points gives
control on the diameters of  $W'_{t-1}$ for those $t$ for which $W_t$ 
contains a critical point satisfying (CE), using (2.1) as below. 
For $\delta$ small enough we also control diameters of $W$'s containing 
critical points not satisfying (CE) (i. e. whose forward trajectory contains other critical points).

Thus by (2.1), for $t>0$ the smallest integer
such that $f^t(c)$ is not critical for iterates of $f,$ 
applied to  $F=f^{s-t}$  (we can assume
$s-t$ positive because it is sufficient to consider only $s$ large)

$$|{f^{(s-t)}}'(f^t(c))|\diam(W'_{s-t}) \leq C_1 b_{s-t+1}^{-C_2} \d.$$


Now (CE) gives, for every $\theta>1,$

$$\dist(c,f^{n-s}(x)) \leq (C_\theta \l^{-s} \theta^s\d) ^\frac1\nu
 < \l^{-\frac{s}{2\nu}}$$
if $\d$ is small enough.
We obtain the contradiction $n\in S_{n-s}$ and conclude that
(with $D=\nu N_f$) $A\subset G(x).$
\karo\enddemo

\demo{Proof of Theorem 1.1} We need to pass from good times for $x\in J$
to good scales, in which $J^c$ contains some definite disc. The argument we
use is similar to the proof of Lemma 1 in \cite{LP}.

Let $\d$ be as in Lemma 2.2 and let $W_n(x)$ be the component of 
$f^{-n}( B(f^n(x),\d/2) )$ containing $x.$ 

Denote by $r(W_n(x))=\dist(x,\partial W_n(x))$ the inradius of $W_n(x).$
We claim that there is an integer $N$ such that the following holds:
For all $x\in J$ and for all $n,n'\in G(x)$ with $n-n'\geq N$
$$\diam(W_n)\leq \frac12 r(W_{n'}).\tag2.5$$

For $n'=0$ ($W_0=B(x,\d/2)$) this is essentially Ma\~n\'e's result \cite{M},
see \cite{P1, Lemma 1.1}. In fact, for each $0<\tau<1$ there is
$N=N(\tau,f,D,\d)$ such that $\diam(W_n(x))\leq \tau\d/2.$

For $n'>0$ use backward iteration: As $f^n$ is $D-$ critical on $W_n(x),$
$f^{(n-n')}$ is $D-$ critical on $W_{n-n'}(f^n(x))$, so that
$$f^{n'}(W_n) = W_{n-n'}(f^{n'}(x))\subset B(f^{n'}(x),\tau\frac\d2)$$
by the first case.
Applying $f^{-n'}$ we obtain (2.5) provided $\tau$ is small enough,
by (2.2).

Consider the increasing sequence $g_j$ of all good times of $x,$
$\{g_j\}=G(x),$ and set $k_j=g_{Nj}.$ By Lemma 2.2 we have
$k_j\leq2 N j,$ and as $k_{j+1}-k_j\geq N$ inequality (2.5) implies
$$\diam(W_{k_{j+1}}(x)) < \frac12 \diam(W_{k_j}(x)).\tag2.6$$
On the other hand, as $f$ is Lipschitz continuous, there is a constant
$L$ such that $\diam(W_n(x)) > 2^{-n L}$ for all $n$ and $x.$
Hence we obtain an increasing sequence $n_j < p_1 j$ (with $p_1<2 L N$)
such that
$$\diam (W_{k_j}(x)) \sim 2^{-n_j}.$$

As $J$ is nowhere dense, there is $\tau>0$ such for every $y\in J$ there
is a disc $U\subset B(y,\d/2)\setminus J$ of radius $\tau\d/2.$

To show porosity of $J$ at $x\in J,$ apply the last statement to $y=f^{k_j}(x)$
(with the sequence $k_j$ constructed above). By (2.3) we find that a
component $V$ of $f^{-k_j}(U)$ in $W_{k_j}(x)$
contains a disc of radius $\geq C_4 2^{-n_j}.$
For the center $z_j$ of this disc we have
$\dist(z_j,x)<\diam(W_{k_j}(x))\sim 2^{-n_j}$ and
$\dist(z_j,J)> C_4 2^{-n_j}.$
We have thus found the desired sequence of discs in $J^c$ and the proof is finished.
\karo\enddemo

\heading3. On H\"older Fatou components\endheading

As a by-product, Lemma 2.2 gives a new approach to the result of Graczyk and Smirnov:

\proclaim{Theorem \cite{GS}} If $f$ is as in Theorem 1.1 and $A$ is a
Fatou component, then $A$ is a H\"older domain.
\endproclaim

See \cite{SS}\cite{JM}\cite{KR}\cite{GS} 
and the references therein for the definition
and results about H\"older domains.

The main estimate in our proof of the above theorem, replacing the second
\ce\ condition established and used in \cite{GS}, is
\proclaim{Proposition 3.1}There exist $0<\xi<1$ and $\d_0>0$ such that for
all $n,$ all $x\in J$ and
for every component $W$ of $f^{-n}(B(f^n(x),\d_0)$
$$\diam W \leq \xi^n.$$
\endproclaim
\demo{Proof} Using \cite{P1,Remark 3.2} we find an integer $N$ and
$\d_0>0$ (smaller than the $\d$ of Lemma 2.2)
with the property that every component of $f^{-m}(B(f^m(y),\d_0)$
has diameter less than
$\d_0/2$ whenever $m\geq N$ and $y\in J$
($m$ does not have to be a good time for
$y$).

As in the proof of Theorem 1.1 let $k_j\in G(x)$ be the $Nj-$th good
time of $x$ ($k_j=g_{Nj}$ with $\{g_j\}=G(x)$).
For $n>N$ and $W$ as above, let $k$ be the largest of the
$k_j$ with $n-k\geq N.$ Then
$$f^k(W)\subset B(f^k(x), \d_0)$$
and from (2.6) we obtain
$$\diam(W) \leq (\frac12)^\frac{n-N}{2N}.$$
\karo\enddemo

>From Proposition 3.1 the H\"older property of invariant Fatou components
(and thus of all components, \cite{GS, Lemma 5.4}) can be concluded as
in \cite{GS, Section 5}: Let $F$ be an attractive (or superattractive)
invariant Fatou component and pick
$\Omega\Subset F$ open, containing all critical points in $F,$ with
$f(\Omega)\Subset\Omega$ and such that $\dist(x,\Omega)<\d_0$ for
all $x\in\partial F.$

Fix a point $z_0\in \Omega.$
Set $n(z)=\min\{n\geq0:f^n(z)\in\Omega\}$
for $z\in F.$ Then the quasihyperbolic distance satisfies
$$\dist_{qh} (z,z_0)\sim n(z)$$
for $z\in F\setminus\Omega$ by \cite{GS, Lemma 5.2}.

For $z\in F$ with $\dist(z,J)<\d_0,$ let $x\in J$ be closest to
$f^{n(z)-1}(z),$ then \newline
$\dist(f^{n(z)-1}(z),x)<\d_0$ and we obtain
$$\dist(z, J) < \xi^{n(z)-1}$$
by Proposition 3.1.
It follows that 
$$n(z)\sim \dist_{qh}(z,z_0)\lesssim \log\frac1{\dist(z,J)},$$
establishing the H\"older property of $F.$

\newpage

\heading 4. Dimension of porous sets \endheading

In this section, we give a simple proof of an estimate of Minkowski
dimension that is sufficient for the proof of Corollary 1.2, without using the
estimates from \cite{KR}.

\

Let $E\subset\R^d$ be a bounded set.
Call $E$ {\it mean porous in all directions} if 
for all $\a>0$  there exist $\b>0,P>0$ such that for all $z\in E$ there 
exists a sequence 
$n_j\le Pj$ such that for all $j$ and for all balls 
$B(z',\a 2^{-n_j})\subset B(z,2^{-n_j}),$  there exists a ball 
$B(z'',\b 2^{-n_j})\subset B(z',\a 2^{-n_j})\setminus E.$
It is easy to see that there are sets which are mean porous but not
mean porous in all directions.

\proclaim{Theorem 4.1} If $f$ satisfies the \cec , has no parabolic
periodic point and if $J(f)\neq\c,$ then $J(f)$ is mean porous in all
directions.
\endproclaim

\demo{Proof} This is a small modification of the proof of Theorem 1.1.
Only the last two paragraphs of the proof need to be changed: The fact 
that $J$ is nowhere dense guarantees, for each $\overline\a>0,$
the existence of $\overline\b>0$ and disks
$B(y'',\overline\b\d/2)\subset B(y',\overline\a\d/2)\setminus J$ 
for all $y\in J$ and
all disks $B(y',\overline\a\d/2)\subset B(y,\d/2).$ 
Applying this to $\overline\a=C\a$ and
taking preimages as in the
the proof of Theorem 1.1 proves the claim, with $\b=C'\overline\b.$ 
\karo\enddemo

Consider a covering $\B_n$ of $\R^d$ by boxes of the form
$[p_1 2^{-n}, (p_1+1)2^{-n}) \times  ...$\newline $ \times [p_d 2^{-n}, (p_d+1)2^{-n})$ 
for all sequences of integers 
$p_1,p_2,...,p_d$. For every $z\in E$ write $Q(z,n)$ for the element of 
$\B_n$ which contains $z$.
Call $E\subset \R^n$ {\it box mean porous} if there exist  
$N,P>0$ and for all $z\in E$ a sequence $n_j\le Pj$ and 
$Q_j\in \B_{n_j+N}$ such that $Q_j\subset Q(z,n_j)\setminus E$.

Of course, {\it mean porosity in all directions} implies 
{\it box mean porosity}, which in turn implies {\it mean porosity}.
Thus the next statement follows also from \cite{KR}.

\proclaim{Proposition 4.2} If a bounded set $E\subset\R^d$ is
box mean porous, then the Minkowski dimension of $E$ satisfies
$\operatorname{MD}(E)<d$.
\endproclaim

\demo{Proof} We may assume $E$ is contained in the unit box $[0,1)^d=:Q_0$. 
Consider the graph (a tree) $\T$ whose vertices are those elements of all 
$\B_n, n=0,1,...,$ which intersect $E$. 
We join $Q\in \B_n$ to $Q'\in\B_{n+1}$ with 
an edge if $Q\supset Q'$. For every vertex $Q\in \T$ call all the vertices 
in the line in $|\T|$ (the body of $\T$) joining $Q$ to $Q_0$ {\it 
ancestors}. If $Q\in\B_n$ and $Q'\in\B_{n+k}$ with $Q'\subset Q,$
let us say that
$Q'$ a $k$-{\it child} of $Q.$

Denote $K=(2^d)^N$. We shall prove that for any integer $n$ being a multiple
of $PN,$ the number of $n$-children of $Q_0$ is at most 
$C (K-1)^{n/(PN)}K^{n/N-n/(PN)},$ where $C$ depends on $d$ and $N$ only. 
As this equals  $ C (2^d)^{\a n}$ with 
$\a=1-{1-\log(K-1)/\log K \over P},$ the estimate $\operatorname{MD}(E)\le \a d$
follows at once.

\
\noindent
The definition of box mean porosity implies two properties of $\T$:

(i)  The number of $N$-children of each vertex is at most $K$; 

(ii) For all $v\in \B_n$ there exists a sequence of at least $n/P$ ancestors
each of which has at most $K-1$  $N$-children.


We want to estimate $\# \B_n$ from above, assuming (i) and (ii). 
The following simple argument is due to Michal Rams.

For every $0\le b< N,$ define a measure $\m_b$ on $\B_n$ by inductively
distributing mass from the top of $\T$ down to $\B_n$ as follows:
Start with equidistributing unit mass on $\B_b$. Given a mass at each vertex 
$v\in \B_{kN+b},$ equidistribute  it on $\Ch_N(v)$, where $\Ch_N(v)$
denotes the set of $N$-children of $v$. 
At the last step, if $n-(kN+b)=k'<N,$ we equidistribute on $\Ch_{k'}(v)$. 
Set $\mu = \frac1N\sum_{b=0}^{N-1} \mu_b$. Then for every $v\in \B_n$ 
there is at least one $b=b(\nu)$ with 
$\mu_b(v) \geq (K-1)^{n/(PN)-N}K^{n/N-n/(PN)+N},$ 
which gives the demanded upper bound on $\# \B_n$.
\karo\enddemo

\Refs
\widestnumber\key{McM}


[BC] M. Benedicks, L. Carleson: The dynamics of the H\'enon maps.
Ann. of Math. 133 (1991), 73-169.

[CE] P. Collet, J.-P. Eckmann: Positive Lyapunov exponents and absolute 
continuity for maps of the interval. Ergodic Th. and Dyn. Sys. 3 (1983), 13-46.

[CJY] L. Carleson, P. Jones, J.-C. Yoccoz: Julia and John. 
Bol. Soc. Bras. Mat. 25 (1994), 1-30.

[DU] M. Denker, M. Urba\'nski: On Hausdorff measures on Julia sets of
subexpanding rational maps. Isr. J. Math. 76 (1992), 193-214.

[DPU] M. Denker, F. Przytycki, M. Urba\'nski: On the transfer operator for 
rational functions on the Riemann sphere. Preprint SFB 170 G\"ottingen, 
4 (1994). To appear in Ergodic Th. and Dyn. Sys.

[G]  J. Graczyk: Hyperbolic subsets, conformal measures and Hausdorff dimension
of Julia sets. Preprint, 1995.

[GS] J. Graczyk, S. Smirnov: Collet, Eckmann, \& H\"older. Preprint, 
November 1995 and February 1996.

[JM] P. Jones, N. Makarov: Density properties of harmonic measure.
Ann. of Math. 142 (1995), 427-455.

[KR] P. Koskela, S. Rohde: Hausdorff dimension and mean porosity. 
Preprint, 1995

[LP] G. Levin, F. Przytycki: When do two rational function have the same
Julia set? Preprint, to appear in Proc. AMS.

[M] R. Ma\~n\'e: On a theorem of Fatou. Bol. Soc. Bras. Mat. 24 (1993), 1-12.

[McM] C. McMullen: Self-similarity of Siegel discs and Hausdorff dimension
of Julia sets. Preprint, 1995.

[NP] T. Nowicki, F. Przytycki: The conjugacy of Collet-Eckmann's map of the
interval with the tent map is H\"older continuous. 
Ergodic Th. and Dyn. Sys. 9 (1989), 379-388.

[NPV] F. Nazarov, I. Popovici, A. Volberg: Domains with quasi-hyperbolic
boundary condition, domains with Greens boundary condition, and estimates of
Hausdorff dimension of their boundary. Preprint, 1995.

[P] Chr. Pommerenke: Boundary behaviour of conformal maps. Springer, 1992.

[P1] F. Przytycki: Iterations of holomorphic Collet-Eckmann maps: conformal 
and invariant measures. Preprint n$^o$ 57, 
Lab. Top. Universit\'e de Bourgogne, F\'evrier 1995.

[P2] F. Przytycki: On measure and Hausdorff dimension of Julia sets 
for holomorphic Collet-Eckmann maps. MSRI Berkeley Preprint 072-95
(July 1995). 

[PUZ] F. Przytycki, M. Urba\'nski, A. Zdunik: Harmonic, Gibbs and Hausdorff
measures for holomorphic maps, Part 2. Studia Math. 97 (1991), 189-225.

[SS] W. Smith, D. Stegenga: Exponential integrability of the quasihyperbolic
metric in H\"older domains. Ann. Acad. Sci. Fenn. 16 (1991), 345-360.

[T] M. Tsujii: Positive Lyapunov exponents in families of one dimensional 
dynamical systems. Invent. Math. 111 (1993), 113-137.

[U] M. Urba\'nski: Rational functions with no recurrent critical points. 
Ergodic Th. and Dyn. Sys. 14.2 (1994), 391-414.

\endRefs

\enddocument